\documentclass[11pt]{article}
\usepackage{amsmath}
\usepackage{amssymb}
\usepackage{enumerate}
\usepackage{color}
\usepackage{xfrac}

\newtheorem{thm}{Theorem}[section]

\newtheorem{lemma}[thm]{Lemma}
\newtheorem{prop}[thm]{Proposition}
\newtheorem{rem}[thm]{Remark}
\newtheorem{cor}[thm]{Corollary}

\newcommand{\cc}[1]{c_{#1}^{\lambda}}
\newcommand{\g}[1]{\gamma_{#1}^{\lambda}}

\begin{document}

\title{On generalized Stieltjes functions}

\author{Stamatis Koumandos and Henrik L. Pedersen\footnote{Research supported by grant DFF–4181-00502 from The Danish Council for Independent Research $|$ Natural Sciences}
}

\date{\today}
\maketitle

\begin{abstract}
It is shown that a function $f$ is a generalized Stieltjes function of order $\lambda>0$ if and only if $x^{1-\lambda}(x^{\lambda-1+k}f(x))^{(k)}$ is completely monotonic for all $k\geq 0$, thereby complementing a result due to Sokal.
Furthermore, a characterization of those completely monotonic functions $f$ for which $x^{1-\lambda}(x^{\lambda-1+k}f(x))^{(k)}$ is completely monotonic for all $k\leq n$ is obtained in terms of properties of the representing measure of $f$.
\end{abstract}


\section{Introduction}

In this paper we investigate a real-variable characterization of generalized Stieltjes functions obtained by Sokal, see \cite{Sok}.

Let $\lambda >0$ be given. A function $f:(0,\infty)\to \mathbb R$ is called a generalized Stieltjes function of order $\lambda$ if 
$$
f(x)=\int_0^{\infty}\frac{d\mu(t)}{(x+t)^{\lambda}}+c,
$$
where $\mu$ is a positive measure on $[0,\infty)$ making the integral converge for $x>0$ and $c\geq 0$. 

The class of ordinary Stieltjes functions is the class of generalized Stieltjes functions of order $1$. 

A $C^{\infty}$-function $f$ 	on $(0,\infty)$ is completely monotonic if $(-1)^nf^{(n)}(x)\geq 0$ for all $n\geq 0$ and all $x>0$. 
Bernstein's theorem characterizes these functions as Laplace transforms of positive measures: $f$ is completely monotonic if and only if there exists a positive measure $\mu$ on $[0,\infty)$ such that $t\mapsto e^{-xt}$ is integrable w.r.t.\ $\mu$ for all $x>0$ and
$$
f(x)=\int_0^{\infty}e^{-xt}\, d\mu(t).
$$

We remark that $f$ is a generalized Stieltjes function of order $\lambda$ if and only if 
\begin{equation}
\label{eq:sokalabsolutely}
f(x)=\int_0^{\infty}e^{-xt}t^{\lambda-1}\varphi(t)\, dt+c, \quad x>0
\end{equation}
for some completely monotonic function $\varphi$, and some non-negative number $c$. See \cite[Lemma 2.1]{kp2}.

Sokal (see \cite{Sok}) introduced for $\lambda >0$ the operators
$$
T_{n,k}^{\lambda}(f)(x)\equiv (-1)^nx^{-(n+\lambda-1)}\left(x^{k+n+\lambda-1}f^{(n)}(x)\right)^{(k)}, \quad n,k\geq 0
$$
and obtained the following characterization.
\begin{thm}
The following are equivalent for a $C^\infty$-function $f$ defined on $(0,\infty)$.
\begin{enumerate}[(a)]
 \item $f$ is a generalized Stieltjes function of order $\lambda$;
\item $T_{n,k}^{\lambda}(f)(x)\geq 0$ 
for all $x>0$, and  $n,k\geq 0$.
\end{enumerate}
  \end{thm}
Sokal's characterization is an extension of Widder's characterization of the class of ordinary Stieltjes functions: $f$ is a Stieltjes function if and only if the function $(x^kf(x))^{(k)}$ is completely monotonic for all $k\geq 0$. 
(See \cite{Wid1}.)

In \cite[Theorem 1.5]{kp3} an analogue of
Sokal's result where the function $\varphi$ in \eqref{eq:sokalabsolutely} is absolutely monotonic is obtained. See also \cite[Theorem 2]{karpprilepkina} for a result complementing \cite[Theorem 1.1]{kp3}. 

\begin{rem}
\label{rem:equivalent}
 Notice that, by Leibniz' rule,
 $$
 x^{-(n+\lambda-1)}\left(x^{k+n+\lambda-1}f^{(n)}(x)\right)^{(k)}=\sum_{j=0}^{k}\binom{k}{j}\frac{\Gamma(n+k+\lambda)}{\Gamma(n+j+\lambda)}x^jf^{(n+j)}(x).
 $$
\end{rem}

In this paper we first show that condition (b) in Sokal's theorem above can be replaced by the condition that 
$$
c_k^{\lambda}(f)(x)\equiv x^{1-\lambda}(x^{\lambda-1+k}f(x))^{(k)}
$$
is completely monotonic for all $k$. There is a simple relation between $T^{\lambda}_{n,k}(f)$ and $\cc{k}(f)$:
\begin{prop} 
\label{prop:direct}
The relation
 $$
 T^{\lambda}_{n,k}(f)(x)=(-1)^n\cc{k}(f)^{(n)}(x)
 $$
holds for any $n,k\geq 0$ and $x>0$.
\end{prop}
\begin{cor}
 \label{cor:sokal} The following are equivalent for a function  $f\in C^{\infty}((0,\infty))$.
 \begin{enumerate}
  \item[(i)] $f$ is a generalized Stieltjes function of order $\lambda$;
  \item[(ii)] $\cc{k}(f)$
is completely monotonic for all $k\geq 0$. 
 \item[(iii)]$T^{\lambda}_{n,k}(f)\geq 0$ for all $n\geq 0$ and all $k\geq 0$.
 \end{enumerate}
\end{cor}

In \cite{kp2} the generalized Stieltjes functions corresponding to measures having moments of all orders were charaterized in terms of properties of remainders in asymptotic expansions. (A measure $\mu$ has moments of all orders in any 
polynomial is integrable w.r.t.\ $\mu$.)
In view of the results in the 
present paper we notice the following corollary. The proof follows by combining Corollary \ref{cor:sokal} with \cite[Theorem 3.1]{kp2} and \cite[Lemma 3.1]{kp2}.
\begin{cor}
 The following are equivalent for a function $f:(0,\infty)\to \mathbb R$.
 \begin{enumerate}
  \item[(i)] $f$ is a generalized Stieltjes function corresponding to a measure $\mu$ having moments of all orders;
  \item[(ii)] $\cc{k}(f)$ is completely monotonic for all $k\geq 0$ and the function $x^{\lambda-1}f(x)$ admits for any $n$ an asymptotic expansion 
  $$
  x^{\lambda-1}f(x)=\sum_{k=0}^{n-1}\frac{\alpha_k}{x^{k+1}}+r_n(x),
  $$
  in which $x^nr_n(x)\to 0$ as $x\to \infty$.
 \end{enumerate}
 In the affirmative case, $\alpha_k=(-1)^k(\lambda)_ks_k(\mu)/k!$ where $s_k(\mu)$ is the $k$'th moment of $\mu$, and $r_n$ has the representation
 $$
 r_n(x)=(-1)^nx^{\lambda-1}\int_0^{\infty}e^{-xt}t^{\lambda-1}\xi_n(t)\, dt,
 $$
 where  $\xi_n$ belongs to $C^{\infty}([0,\infty))$,  and satisfies $\xi_n^{(j)}(0)=0$ for $j\leq n-1$ and $0\leq \xi_n^{(n)}(t)\leq s_n(\mu)$ for $t\geq 0$. Furthermore,
 $$
 \cc{n}(f)(x)=x^{1-\lambda}(x^nr_n(x))^{(n)}=\cc{n}\left(\mathcal L(t^{\lambda-1}(-1)^n\xi_n(t)\right)(x).
 $$
\end{cor}

Our aim is also to characterize, for any given positive integer $N$, those functions $f$ for which $c_0^{\lambda}(f),\ldots,c_N^{\lambda}(f)$ are completely monotonic. % for any given number $N$.
In the case where $\lambda=1$ this has been carried out in \cite{p}, but the case of general $\lambda$ requires, as we shall see, additional insight.

We thus introduce the classes $\mathcal C_N^{\lambda}$ as 
$$
\mathcal C_N^{\lambda}=\{f\in C^{\infty}((0,\infty))\, |\, c_k^{\lambda}(f)\ \text{is completely monotonic for }\ k=0,\ldots,N\}.
$$
We shall use some distribution theory so we briefly describe our notation. The action of a distribution $u$ on a test function $\varphi$ (an infinitely often differentiable function of compact support in $(0,\infty)$) is denoted by $\langle u,\varphi\rangle$. 
The distribution $\partial u$ is defined via $\langle \partial u,\varphi\rangle=-\langle u,\varphi' \rangle$. A standard reference to distribution theory is \cite{rudin}.

Our results can be formulated as follows.

\begin{thm}
\label{thm:main}
Let $\lambda >0$ be given, and let $N\geq 1$. The following properties of a function $f:(0,\infty)\to \mathbb R$ are equivalent.
\begin{enumerate}[(a)]
\item $f\in \mathcal C_N^{\lambda}$;
\item $f$ can be represented as 
$$f(x)=c+\int_0^{\infty}e^{-xs}s^{\lambda-1}\, d\mu(s),$$
where $c\geq 0$, and $\mu$ is a positive measure on $(0,\infty)$ for which $\mu_k\equiv (-1)^ks^k\partial^k\mu$, (in distributional sense)  is a positive measure such that  
\begin{align*}
&\int_0^{\infty}e^{-xs}s^{\lambda-1}\, d\mu_k(s)<\infty, \quad k=0,\ldots,N.
\end{align*}
\end{enumerate}
In the affirmative case,
$$
\cc{k}(f)(x)=x^{1-\lambda}\left(x^{\lambda-1+k}f(x)\right)^{(k)}=\int_0^{\infty}e^{-xs}s^{\lambda-1}\, d\mu_k(s)
$$
for $k=0,\ldots,N$.
\end{thm}

We notice the following corollary characterizing those non-negative functions $f$ for which $\cc{1}(f)$ is completely monotonic. The proof follows from Propostion \ref{prop:specialcase} and Lemma \ref{lemma:integrability}.
\begin{cor}
\label{cor:bernstein}
 Let $f$ be a non-negative $C^{\infty}$-function defined on $(0,\infty)$. Then $x^{1-\lambda}\left(x^\lambda f(x)\right)'=\lambda f(x)+xf'(x)$ is completely monotonic if and only if
 $$
 f(x)=\alpha+\frac{\beta}{x^{\lambda}}+\int_0^{\infty}e^{-xs}s^{\lambda-1}\int_s^{\infty}\frac{d\mu(t)}{t^{\lambda}}\, ds,
 $$
 for some non-negative numbers $\alpha$ and $\beta$ and some positive measure $\mu$ on $(0,\infty)$ making the integral convergent.
\end{cor}

\begin{rem} 
 It is easy to see that 
 $e^{-xs}s^{\lambda-1}\int_s^{\infty}d\mu(t)/t^{\lambda}$ is integrable on $(0,\infty)$ if and only if $s^{\lambda-1}\int_s^{\infty}d\mu(t)/t^{\lambda}$ is integrable at $0$, and that this is the case if and only if 
 $\int_0^1\, d\mu(t)<\infty$ and $\int_1^{\infty}\, d\mu(t)/t^{\lambda}<\infty$.
\end{rem}
Corollary \ref{cor:bernstein} can be reformulated as follows. Let $g$ be a non-negative $C^{\infty}$-function on $(0,\infty)$. Then  $x^{1-\lambda}g'(x)$ is completely monotonic if and only if 
 \begin{equation}
 \label{eq:bernstein}
 g(x)=\alpha x^{\lambda}+\beta+\int_0^{\infty}\int_0^{xt}e^{-u}u^{\lambda-1}\, du\frac{d\mu(t)}{t^{\lambda}}.
 \end{equation}
 Formulated in this way the corollary is related to the class of Bernstein functions.
A Bernstein function is by definition a non-negative function $g$ on $(0,\infty)$ for which $g'$ is completely monotonic. These functions admit an integral representation 
(see \cite[Theorem 3.2]{SSV} or \cite{bergforst}), which we for the reader's convenience state here: $g$ is a Bernstein 
function if and only if 
$$
g(x)=\alpha x+\beta+\int_0^{\infty}(1-e^{-xt})\, d\nu(t),
$$
where $\alpha$ and $\beta$ are non-negative numbers, and $\nu$, called the L\'evy measure, is a positive measure on $(0,\infty)$ satisfying
$\int_{0}^{1}t d\nu(t)<\infty$  and $\int_{1}^{\infty} d\nu(t)<\infty$.

When $\lambda=1$, we have 
 $$
 \int_0^{xt}e^{-u}u^{\lambda-1}\, du=1-e^{-xt},
 $$
 and \eqref{eq:bernstein} reduces to the integral representation of a Bernstein function with the corresponding L\'evy measure being $d\mu(t)/t$. Corollary \ref{cor:bernstein} contains a characterization of what could be called ``generalized Bernstein functions of order $\lambda$''.

\section{Proofs}
{\it Proof of Proposition \ref{prop:direct}:} 
The key to the proof is the following relation
\begin{equation}
\label{eq:key}
 T^{\lambda}_{n,k}(f)(x)=(-1)^n\left(\sum_{j=0}^k(\lambda-1)_{k-j}\binom{k}{j}\left(x^{j}f(x)\right)^{(j)}\right)^{(n)},
\end{equation}
which we verify now. A standard application of Leibniz' formula yields
\begin{align*}
\left(x^{j}f(x)\right)^{(j+n)}&=\sum_{l=0}^{n+j}\binom{n+j}{l}\left( x^j\right)^{(l)}f^{(n+j-l)}(x)\\
&=\sum_{l=0}^{j}\binom{n+j}{l}\frac{j!}{(j-l)!}x^{j-l}f^{(n+j-l)}(x)\\
&=\sum_{m=0}^{j}\binom{n+j}{j-m}\frac{j!}{m!}x^{m}f^{(n+m)}(x).
\end{align*}
Hence, the right hand side of \eqref{eq:key} equals
\begin{align*}
\lefteqn{ (-1)^n\sum_{j=0}^k(\lambda-1)_{k-j}\binom{k}{j}\left(x^{j}f(x)\right)^{(j+n)}}\\
&= (-1)^n\sum_{j=0}^k(\lambda-1)_{k-j}\binom{k}{j}\sum_{m=0}^{j}\binom{n+j}{j-m}\frac{j!}{m!}x^{m}f^{(n+m)}(x)\\
&= (-1)^n\sum_{m=0}^k\left\{\sum_{j=m}^{k}(\lambda-1)_{k-j}\binom{k}{j}\binom{n+j}{j-m}\frac{j!}{m!}\right\}x^{m}f^{(n+m)}(x).
\end{align*}
The expression in the brackets can be written in another form. Indeed 
$$
\sum_{j=m}^{k}(\lambda-1)_{k-j}\binom{k}{j}\binom{n+j}{j-m}\frac{j!}{m!}=\binom{k}{m}\frac{\Gamma(n+k+\lambda)}{\Gamma(n+m+\lambda)},
$$
by a corollary to the Chu-Vandermonde identity (see \cite[p.\ 70]{aar}).
This gives us
\begin{align*}
\lefteqn{(-1)^n\left(\sum_{j=0}^k(\lambda-1)_{k-j}\binom{k}{j}\left(x^{j}f(x)\right)^{(j)}\right)^{(n)}}\\
&=(-1)^n\sum_{m=0}^k\binom{k}{m}\frac{\Gamma(n+k+\lambda)}{\Gamma(n+m+\lambda)}x^{m}f^{(n+m)}(x)=T^{\lambda}_{n,k}(f)(x).
\end{align*}
For $n=0$ the identity reads
$$
\sum_{j=0}^k(\lambda-1)_{k-j}\binom{k}{j}\left(x^{j}f(x)\right)^{(j)}=T^{\lambda}_{0,k}(f)(x)=\cc{k}(f)(x),
$$
and the proposition is proved.\hfill $\square$

To prove Theorem \ref{thm:main} we need a few preliminary results.
\begin{lemma}
\label{lemma:recursion}
For $f\in C^{\infty	}((0,\infty))$ we have
 $$
 \cc{k+1}(f)(x)=(\lambda+k)\cc{k}(f)(x)+x\cc{k}(f)'(x).
 $$
\end{lemma}
{\it Proof.} This follows by computation:
\begin{align*}
 \cc{k+1}(f)(x)&=x^{1-\lambda}(xx^{\lambda-1+k}f(x))^{(k+1)}\\
 &=x^{2-\lambda}(x^{\lambda-1+k}f(x))^{(k+1)}+(k+1)x^{1-\lambda}(x^{\lambda-1+k}f(x))^{(k)}\\
 %&=x^{2-\lambda}(x^{1-\lambda}(x^{\lambda-1+k}f(x))^{(k)}x^{\lambda-1})'+(k+1)\cc{k}(f)(x)\\
 &=x^{2-\lambda}(x^{\lambda-1}\cc{k}(f)(x))'+(k+1)\cc{k}(f)(x)\\
 &=x\cc{k}(f)'(x)+(\lambda+k)\cc{k}(f)(x).
 \end{align*}

 \hfill $\square$
 
 \begin{prop}
 \label{prop:recursion}
  Suppose that $f\in \mathcal C_N^{\lambda}$, and let for $k=0,\ldots,N$
  $$\cc{k}(f)(x)=\int_0^{\infty}e^{-xs}\, d\mu_k(s)+b_k,$$
  where $\mu_k$ is a positive measure on $(0,\infty)$ and $b_k\geq 0$. Then, in the distributional sense,
  $$
  (-1)^ks^k\partial^k(s^{1-\lambda}\mu_0)=s^{1-\lambda}\mu_k.
  $$
 \end{prop}
{\it Proof.} From Lemma \ref{lemma:recursion} it follows that (for $k\leq N-1$)
\begin{align*}
\int_0^{\infty}e^{-xs}\, d\mu_{k+1}(s)+b_{k+1}&=(\lambda+k)\int_0^{\infty}e^{-xs}\, d\mu_k(s)+(\lambda+k)b_k\\
&{}\qquad -x\int_0^{\infty}se^{-xs}\, d\mu_k(s).
\end{align*}
Letting $x\to \infty$ yields $b_{k+1}=(\lambda+k)b_k$ so that
\begin{align*}
\frac{1}{x}\int_0^{\infty}e^{-xs}\, d\mu_{k+1}(s)&=(\lambda+k)\frac{1}{x}\int_0^{\infty}e^{-xs}\, d\mu_k(s)\\
&{}\qquad -\int_0^{\infty}se^{-xs}\, d\mu_k(s).
\end{align*}
By the uniqueness of the Laplace transform we obtain
$$
s\mu_k=((\lambda+k)\mu_k-\mu_{k+1})\ast m.
$$
(Here, $m$ denotes Lebesgue measure on $(0,\infty)$.) We get by differentiation (as distributions) that 
$$
s\partial \mu_k=(\lambda+k-1)\mu_k-\mu_{k+1}.
$$
We shall obtain the assertion in the proposition by induction, using this recursive relation: for $k=0$ the assertion is valid. Before verifying the induction step notice that 
$$
s\partial (s^k\partial^k(s^{1-\lambda}\mu_0))=ks^k\partial^k(s^{1-\lambda}\mu_0)+s^{k+1}\partial^{k+1}(s^{1-\lambda}\mu_0).$$
Suppose now that the assertion holds for $k$. Then
\begin{align*}
 s^{k+1}\partial^{k+1}(s^{1-\lambda}\mu_0)&=s\partial (s^k\partial^k(s^{1-\lambda}\mu_0))-ks^k\partial^k(s^{1-\lambda}\mu_0)\\
 &=s\partial ((-1)^ks^{1-\lambda}\mu_k)-k(-1)^ks^{1-\lambda}\mu_k\\
 &=(-1)^k\{s(1-\lambda)s^{-\lambda}\mu_k+s^{1-\lambda}s\partial \mu_k-ks^{1-\lambda}\mu_k\}\\
 &=(-1)^ks^{1-\lambda}\{(1-\lambda)\mu_k+(\lambda+k-1)\mu_k-\mu_{k+1}-k\mu_k\}\\
 &=(-1)^{k+1}s^{1-\lambda}\mu_{k+1}.
\end{align*}
The assertion holds also for $k+1$, and the proof follows.\hfill $\square$

{\it Proof that (a) implies (b) in Theorem \ref{thm:main}.} If $f\in \mathcal C_N^{\lambda}$ then the function $\cc{k}(f)$ is completely monotonic for $0\leq k\leq N$. In particular
$$
f(x)=\cc{0}(f)(x)=\int_0^{\infty}e^{-xs}\,d\mu_0(s)+b_0=\int_0^{\infty}e^{-xs}s^{\lambda-1}\,d(s^{1-\lambda}\mu_0)(s)+b_0.
$$
Let $\mu=s^{1-\lambda}\mu_0$ and notice that by Proposition \ref{prop:recursion} $(-1)^ks^k\partial^k\mu\, (=s^{1-\lambda}\mu_k)$ is a positive measure with the property that 
$$
\int_0^{\infty}e^{-xs}s^{\lambda-1}\, d((-1)^ks^k\partial^k\mu)(s)<\infty.
$$
Thus (b) follows. \hfill $\square$

The next result is a special case of (b) implies (a) in Theorem \ref{thm:main}. We state and prove it separately in order to describe the method, which will be alluded to in the following proof.   
\begin{prop}
\label{prop:specialcase}
Let $f$ have the representation 
$$
f(x)=c+\int_0^{\infty}e^{-xs}s^{\lambda-1}\, d\mu(s),
$$
where $c\geq 0$ and $\mu$ is a positive measure on $(0,\infty)$. If $-s\partial \mu(s)$ is a positive measure then 
$$
\cc{1}(f)(x)=\lambda c+\int_0^{\infty}e^{-xs}s^{\lambda-1}\, d(-s\partial \mu(s))
$$
is completely monotonic.
\end{prop}
{\it Proof.} Let $n\geq 1$ and take $\varphi_n\in C^{\infty}((0,\infty))$ such that $0\leq \varphi_n(t)\leq 1$,
$$
\varphi_n(t)=\left\{\begin{array}{ll}
             0,&t<1/(2n)\\
             1,&1/n\leq t\leq n\\
             0,&t>2n
            \end{array}\right.,
$$
$|\varphi_n'(t)|\leq \text{Const}\, n$ for $t\in(1/(2n),1/n)$, and $|\varphi_n'(t)|\leq \text{Const}$ for $t\in(n,2n)$.
By definition of the derivative in distributional sense we have 
\begin{align*}
 \lefteqn{\int_0^{\infty}e^{-xs}s^{\lambda-1}\varphi_n(s)\, d(-s\partial \mu(s))}\\
 &=\langle -s\partial \mu(s),e^{-xs}s^{\lambda-1}\varphi_n(s)\rangle\\
 &=\langle \mu,(e^{-xs}s^{\lambda}\varphi_n(s))'\rangle\\
 &=-x\int_0^{\infty}e^{-xs}s^{\lambda}\varphi_n(s)\, d\mu(s)+\int_0^{\infty}e^{-xs}\lambda s^{\lambda-1}\varphi_n(s)\, d\mu(s)\\
 &\quad +\int_0^{\infty}e^{-xs}s^{\lambda}\varphi_n'(s)\, d\mu(s).
 \end{align*}

Using dominated convergence it follows that the sum of first and second term on the right hand side tends to 
$$
-x\int_0^{\infty}e^{-xs}s^{\lambda}\, d\mu(s) +\lambda \int_0^{\infty}e^{-xs}s^{\lambda-1}\, d\mu(s).
$$
The third term tends to zero, again due to dominated convergence and the estimate (using $|s\varphi'(s)|\leq \text{Const}$ for $s\leq 1/n$)
\begin{align*}
\lefteqn{\int_0^{\infty}\left|e^{-xs}s^{\lambda}\varphi_n'(s)\right|\, d\mu(s)\leq}\\ & \text{Const}\left(\,\int_{\sfrac{1}{2n}}^{\sfrac{1}{n}}e^{-xs}s^{\lambda-1}\, d\mu(s)+\int_{n}^{2n}e^{-xs}s^{\lambda}\, d\mu(s)\right).
\end{align*}
Hence, letting $n$ tend to infinity, we obtain that
\begin{align*}
 \cc{1}(f)(x)&=\lambda f(x)+xf'(x)\\
 &=\lambda c+\lambda\int_0^{\infty}e^{-xs}s^{\lambda-1}\, d\mu(s)-x\int_0^{\infty}e^{-xs}s^{\lambda}\, d\mu(s)\\
 &=\lambda c+\int_0^{\infty}e^{-xs}s^{\lambda-1}\, d(-s\partial \mu(s)).
 \end{align*}
Thus $\cc{1}(f)$ is completely monotonic, and $e^{-xs}s^{\lambda-1}$ is integrable w.r.t.\ the measure $-s\partial \mu(s)$.\hfill $\square$

{\it Proof that (b) implies (a) in Theorem \ref{thm:main}.} We suppose that $f$ has the representation 
$$
f(x)=c+\int_0^{\infty}e^{-xs}s^{\lambda-1}\, d\mu(s),
$$
with $c\geq 0$, and $\mu_k\equiv (-1)^ks^k\partial^k \mu$ being a positive measure for $k=0,\ldots,N$. It is easy to verify  that $\mu_{k+1}=k\mu_k-s\partial \mu_k$ for $k=0,\ldots,N-1$. Proposition \ref{prop:specialcase} yields that $\cc{1}(f)$ is 
completely monotonic and has the representation 
$$
\cc{1}(f)(x)=\lambda c+\int_0^{\infty}e^{-xs}s^{\lambda-1}\, d\mu_1(s).
$$
Let us now assume that $\cc{k}(f)$ is completely monotonic
and has the representation 
$$
\cc{k}(f)(x)=b+\int_0^{\infty}e^{-xs}s^{\lambda-1}\, d\mu_k(s).
$$
Then 
\begin{align*}
\cc{k+1}(f)(x)&=(\lambda+k)\cc{k}(f)(x)+x\cc{k}(f)'(x)\\
&=(\lambda+k)\int_0^{\infty}e^{-xs}s^{\lambda-1}\, d\mu_k(s)-x\int_0^{\infty}e^{-xs}s^{\lambda}\, d\mu_k(s).
\end{align*}
Now, taking $\varphi_n$ as before it follows that
\begin{align*}
 \lefteqn{\int_0^{\infty}e^{-xs}s^{\lambda-1}\varphi_n(s)\, d\mu_{k+1}(s)}\\
 &=\int_0^{\infty}e^{-xs}s^{\lambda-1}\varphi_n(s)\, d(k\mu_k-s\partial \mu_k(s))\\
 &=\langle k\mu_k-s\partial \mu_k(s),e^{-xs}s^{\lambda-1}\varphi_n(s)\rangle\\
 &=\langle k\mu_k,e^{-xs}s^{\lambda-1}\varphi_n(s)\rangle+\langle \mu_k,(e^{-xs}s^{\lambda}\varphi_n(s))'\rangle\\
 &=-x\int_0^{\infty}e^{-xs}s^{\lambda}\varphi_n(s)\, d\mu_k(s)+(k+\lambda)\int_0^{\infty}e^{-xs} s^{\lambda-1}\varphi_n(s)\, d\mu_k(s)\\
 &\quad -\int_0^{\infty}e^{-xs}s^{\lambda}\varphi_n'(s)\, d\mu_k(s).
 \end{align*}
 As before, letting $n$ tend to infinity, and applying dominated convergence we get that
$$
\cc{k+1}(f)(x)=\int_0^{\infty}e^{-xs}s^{\lambda-1}\, d\mu_{k+1}(s)
$$
is completely monotonic.\hfill $\square$

\section{Additional results and comments}
Suppose that $\cc{k}(f)$ is completely monotonic for some $k\geq 1$. What can be said about the functions $\cc{0}(f),\ldots,\cc{k-1}(f)$? Are they also completely monotonic? The answer is given in Proposition \ref{prop:question}.

\begin{prop}
\label{prop:question}
Let $k\geq 1$, $f\in C^{\infty}((0,\infty))$ and suppose that the functions $\cc{0}(f),\ldots,\cc{k-1}(f)$ are non-negative. If $\cc{k}(f)$ is completely monotonic then $\cc{j}(f)$ is also completely monotonic for $j\leq k-1$. 
\end{prop}
The proof of this proposition requires some preliminary results. Define 
$$
\g{k}(f)(x)\equiv x^{-\lambda}(x^{\lambda-1+k}f(x))^{(k-1)},\quad k\geq 1.
$$
Notice that $\g{k}(f)=c^{\lambda +1}_{k-1}(f)$. 
\begin{lemma}
\label{lemma:first}
 For $f\in C^{(k+1)}((0,\infty))$ we have
 $$
 \g{k}(f)(x)=\cc{k-1}(f)(x)+(k-1)\g{k-1}(f)(x).
 $$
\end{lemma}
{\it Proof:} This follows by a direct computation.\hfill $\square$

The next lemma is an immediate consequence of Lemma \ref{lemma:first}. 
\begin{lemma}
If $\cc{j}(f)(x)\geq 0$ for all $j=0,\ldots,k$ then $\g{j}(f)(x)\geq 0$ for all $j=1,\ldots,k+1$
\end{lemma}

\begin{lemma}
\label{lemma:comp}
 Let $k\geq 1$ be given and assume that $\cc{j}(f)(x)\geq 0$ for all $j=0,\ldots,k$. Then:
 \begin{enumerate}
  \item[(i)] $(x^{\lambda-1+k}f(x))^{(j)}\geq 0$ for all $j=0,\ldots,k$;
  \item[(ii)] $\lim_{x\to 0}(x^{\lambda-1+k}f(x))^{(j)}= 0$ for all $j=0,\ldots,k-2$;
  \item[(iii)] $\lim_{x\to 0}(x^{\lambda-1+k}f(x))^{(k-1)}\in [0,\infty)$;
 \end{enumerate}

\end{lemma}
{\it Proof.} We use induction in $k$. For $k=1$ (i) is clearly satisfied, (ii) needs not be checked, and (iii) follows by noticing that $x^\lambda f(x)$ is non-negative and increasing. 
For $k=2$, $(x^{\lambda+1}f(x))''=x^{\lambda-1}\cc{2}(f)(x)\geq 0$, $(x^{\lambda+1}f(x))'=x^{\lambda}\cc{1}(f)(x)+\lambda x^\lambda \cc{0}(f)(x)\geq 0$, and thus (i) is satisfied. Property (ii) is clearly satisfied, and (iii) follows since 
$(x^{\lambda+1}f(x))'$ is non-negative and increasing.

Next we assume that $f$ satisfies $\cc{j}(f)\geq 0$ for all $j\leq k+1$, and aim at verifying (i), (ii), and (iii) with $k$ replaced by $k+1$. For $j=k+1$ we  get $(x^{\lambda+k}f(x))^{(j)}=x^{\lambda-1}\cc{k+1}(f)(x)\geq 0$.
For $1\leq j\leq k$ we use
$$
(x^{\lambda+k}f(x))^{(j)}=x(x^{\lambda-1+k}f(x))^{(j)}+j(x^{\lambda-1+k}f(x))^{(j-1)}\geq 0,
$$
and (i) is verified. To see (ii), notice that  
$$
(x^{\lambda+k}f(x))^{(k-1)}=x(x^{\lambda-1+k}f(x))^{(k-1)}+(k-1)(x^{\lambda-1+k}f(x))^{(k-2)},
$$
The last term tends to zero by the induction hypothesis, and the first term equals $x$ times a non-negative and increasing function. Hence (ii) holds for $k+1$. Property (iii) for $k+1$ follows since $(x^{\lambda+k}f(x))^{(k)}$ is a 
positive and increasing function. This proves the lemma.\hfill $\square$
\begin{lemma}
\label{lemma:integrability}
Let $f\in C^{\infty}((0,\infty))$ and suppose that $\cc{0}(f),\ldots,\cc{k-1}(f)$ are non-negative functions. If $\cc{k}(f)$ is completely monotonic then $\g{k}(f)$ is also completely monotonic and 
$$
\g{k}(f)(x)=\frac{l_k}{x^{\lambda}}+\frac{b_k}{\lambda}+\int_0^{\infty}M_k(u)u^{\lambda-1}e^{-xu}\, du,
$$
where
$$
M_k(u)=\int_u^{\infty} s^{-\lambda}d\mu_k(s).
$$

\end{lemma}
{\it Proof.} By the complete monotonicity we may write
$$
(x^{\lambda-1+k}f(x))^{(k)}=x^{\lambda-1}\cc{k}(f)(x)=x^{\lambda-1}\int_0^{\infty}e^{-xs}\, d\mu_k(s)+b_kx^{\lambda-1},
$$
where $b_k\geq 0$ and $\mu_k$ is a positive measure on $(0,\infty)$. The assumptions on non-negativity yield that the function $x^{\lambda}\g{k}(f)(x)=(x^{\lambda-1+k}f(x))^{(k-1)}$ is non-negative and increasing. Hence
$$
l_k\equiv \lim_{x\to 0}x^{\lambda}\g{k}(f)(x)\geq 0.
$$
Furthermore,
\begin{align*}
 x^{\lambda}\g{k}(f)(x)-l_k&=\int_0^x(t^{\lambda-1+k}f(t))^{(k)}\,dt\\
&=\int_0^xt^{\lambda-1}\left(\int_0^{\infty}e^{-ts}\, d\mu_k(s)+b_k\right)\,dt\\
&=\frac{b_k}{\lambda}x^{\lambda}+\int_0^{\infty}\int_0^xt^{\lambda-1}e^{-ts}\, dt\, d\mu_k(s)\\
&=\frac{b_k}{\lambda}x^{\lambda}+x^{\lambda}\int_0^{\infty}\int_0^su^{\lambda-1}e^{-xu}\, du\, s^{-\lambda}d\mu_k(s)\\
&=\frac{b_k}{\lambda}x^{\lambda}+x^{\lambda}\int_0^{\infty}\int_u^{\infty} s^{-\lambda}d\mu_k(s)u^{\lambda-1}e^{-xu}\, du,
\end{align*}
by Fubini's theorem. Consequently,
$$
M_k(u)=\int_u^{\infty} s^{-\lambda}d\mu_k(s)
$$
is finite and $M_k(u)u^{\lambda-1}$ is integrable at $0$. 

The formulas above also show that 
$$
\g{k}(f)(x)=\frac{l_k}{x^{\lambda}}+\frac{b_k}{\lambda}+\int_0^{\infty}M_k(u)u^{\lambda-1}e^{-xu}\, du
$$
is completely monotonic.\hfill $\square$

{\it Proof of Proposition \ref{prop:question}:} From Lemma \ref{lemma:integrability},
$$
(x^{\lambda-1+k}f(x))^{k-1}=l_k+\frac{b_k}{\lambda}x^{\lambda}+x^{\lambda}\int_0^{\infty}M_k(u)u^{\lambda-1}e^{-xu}\, du,
$$
where $l_k,b_k\geq 0$ and $M_k(u)=\int_u^{\infty}s^{-\lambda}\, d\mu(s)$. Notice that 
\begin{equation}
\label{eq:limk}
 u^{\lambda}M_k(u)\to 0, \quad u\to 0.
\end{equation}
(To see this, rewrite as follows 
$$
u^{\lambda}M_k(u)=\int_u^1\left(\frac{u}{s}\right)^{\lambda}\, d\mu(s)+u^{\lambda}\int_1^{\infty}\frac{d\mu(s)}{s^{\lambda}},%\to 0, \quad u\to 0.
$$
and use the dominated convergence theorem on the first term.)

Integrating this relation from $\epsilon$ to $x$, and letting $\epsilon$ tend to 0 we get, using (ii) of Lemma \ref{lemma:comp}, that
\begin{align*}
(x^{\lambda-1+k}f(x))^{(k-2)}&=l_kx+\frac{b_k}{\lambda(\lambda+1)}x^{\lambda+1}\\
&{}\qquad +x^{\lambda+1}\int_0^{\infty}M_{k-1}(u)u^{\lambda}e^{-xu}\, du,
\end{align*}
where 
$$
M_{k-1}(u)=\int_u^{\infty}\frac{M_k(s)}{s^2}\, ds.
$$
Continuing this process (using in each step (ii) of Lemma \ref{lemma:comp}) we get
\begin{align}
x^{\lambda-1+k}f(x)&=\frac{l_k}{(k-1)!}x^{k-1}+\frac{b_k}{(\lambda)_k}x^{\lambda+k-1}\nonumber\\
&{}\qquad +x^{\lambda+k-1}\int_0^{\infty}M_1(u)u^{\lambda+k-2}e^{-xu}\, du, \label{eq:above1}
\end{align}
where 
$$
M_{j}(u)=\int_u^{\infty}\frac{M_{j+1}(s)}{s^2}\, ds, \quad j=1,\ldots,k-1.
$$
Division by $x^{\lambda-1+k}$ in \eqref{eq:above1} shows that $f$ is completely monotonic, and has the representation
\begin{align}
f(x)&=\frac{l_k}{(k-1)!}x^{-\lambda}+\frac{b_k}{(\lambda)_k}+\int_0^{\infty}M_1(u)u^{\lambda+k-2}e^{-xu}\, du\nonumber \\
&=\frac{b_k}{(\lambda)_k}+\int_0^{\infty}\left(M_1(u)u^{k-1}+\frac{l_k}{(k-1)!}\right)e^{-xu}u^{\lambda-1}\, du.\label{eq:conclusionquestion}
\end{align}
In order to show that the functions $\cc{1}(f),\ldots,\cc{k-1}(f)$ are completely monotonic it suffices (Theorem \ref{thm:main}) to verify that
$(-1)^j\partial^j \left(M_1(u)u^{k-1}\right)\geq 0$ for $j=1,\ldots,k-1$. Now, 
$$
M_1(u)u^{k-1}=u^{k-1}\int_u^{\infty}\frac{M_2(s)}{s^2}\, ds=\int_1^{\infty}\frac{M_2(ut)(ut)^{k-2}}{t^k}\, dt,
$$
so it is enough to verify that $(-1)^j\partial^j\left(M_2(u)u^{k-2}\right)\geq 0$ in order to obtain that $(-1)^j\partial^j \left(M_1(u)u^{k-1}\right)\geq 0$. Repeating this argument we end up having to verify that 
\begin{equation}
\label{eq:star}
(-1)^{j}\partial^{j} \left(M_{k-j+1}(u)u^{j-1}\right)\geq 0,\quad 1\leq j\leq k-1.
\end{equation}
These inequalities are verified using induction. For $j=1$ it reads $\partial M_k(u)\leq$ which is true since $M_k$ is a decreasing function. Next assuming that \eqref{eq:star} holds for some $j\leq k-2$ we aim at verifying it for $j+1$. 
We rewrite
the expression $(-1)^{j+1}\partial^{j+1} \left(M_{k-j}(u)u^{j}\right)$ in two ways:
\begin{align*}
 (-1)^{j+1}\partial^{j+1} \left(M_{k-j}(u)u^{j}\right)&=(-1)^{j+1}u\partial^{j+1} \left(M_{k-j}(u)u^{j}\right)\\
 &{}\quad +(-1)^{j+1}(j+1)\partial^{j} \left(M_{k-j}(u)u^{j}\right);\\
 (-1)^{j+1}\partial^{j+1} \left(M_{k-j}(u)u^{j}\right)&=(-1)^{j+1}\partial^{j} \left(-\frac{M_{k-j+1}(u)}{u^2}u^{j+1}\right)\\
 &{}\quad +(-1)^{j+1}(j+1)\partial^{j} \left(M_{k-j}(u)u^{j}\right).
\end{align*}
Comparing these two identities we infer that
$$
(-1)^{j+1}u\partial^{j+1} \left(M_{k-j}(u)u^{j}\right)=(-1)^{j}\partial^{j} \left(M_{k-j+1}(u)u^{j-1}\right),
$$
and thus \eqref{eq:star} holds for $j+1$.\hfill $\square$

\begin{rem} 
Introducing the functions $N_j(u)\equiv M_j(1/u)$ for $j=1,\ldots,k$ it follows that  
 $$
 N_k(u)=\int_{1/u}^{\infty}s^{-\lambda}\, d\mu(s)=\int_{0}^{u}t^{\lambda}\, d\widehat{\mu}(t),
 $$
 where $\widehat{\mu}$ denotes the image measure $\phi(\mu)$, with $\phi(x)=1/x$. For $j\leq k-1$ the relation between $N_j$ and $N_{j+1}$ is
 $$
 N_j(u)=\int_{1/u}^{\infty}\frac{M_{j+1}(s)}{s^2}\, ds=\int_{0}^{u}N_{j+1}(t)\, dt.
 $$
 Consequently we see that the derivatives $N_1^{(j)}(u)$ for $j\leq k-1$ are all non-negative, and take the value $0$ at $u=0$.
 In terms of these functions the representation \eqref{eq:conclusionquestion} can be rewritten as
 \begin{align*}
 f(x)&=\frac{b_k}{(\lambda)_k}+\frac{l_k\Gamma(\lambda)}{(k-1)!x^{\lambda}}+\int_0^{\infty}M_1(u)u^{k-1}e^{-xu}u^{\lambda-1}\, du\\
 &=\frac{b_k}{(\lambda)_k}+\frac{l_k\Gamma(\lambda)}{(k-1)!x^{\lambda}}+\int_0^{\infty}N_1(s)s^{-k-\lambda}e^{-x/s}\, ds.
 \end{align*}
 
\end{rem}

The next proposition shows that for any given $N$ the classes $\mathcal C_N^{\lambda}$ become larger as $\lambda$ increases. As remarked in \cite{Sok} it is not clear how to verify this even for $N=\infty$ only considering the operators $T_{n,k}^{\lambda}$. 
\begin{prop}
If $\lambda_1<\lambda_2$ then $\mathcal C_{N}^{\lambda_1}\subset \mathcal C_{N}^{\lambda_2}$ for all $N$. 
\end{prop}
{\it Proof.} This follows from Leibniz' formula. Assume $f\in \mathcal C_N^{\lambda_1}$ and let $\lambda_2>\lambda_1$. Then 
$$
f(x)=c+\int_0^{\infty}e^{-xs}s^{\lambda_1 -1}\, d\mu(s)=c+\int_0^{\infty}e^{-xs}s^{\lambda_2 -1}\, s^{\lambda_1-\lambda_2}d\mu(s),
$$
where $(-1)^js^{j}\partial^j\mu\geq 0$ for all $j\leq N$. Hence, for $k\leq N$, 
\begin{align*}
 &(-1)^ks^k\partial^k(s^{\lambda_1-\lambda_2}\mu)\\
 &=(-1)^ks^k\sum_{j=0}^k\binom{k}{j}(-1)^{k-j}\frac{\Gamma(k-j+\lambda_2-\lambda_1)}{\Gamma(\lambda_2-\lambda_1)}s^{\lambda_1-\lambda_2+j-k}\partial^j\mu\\
 &=s^{\lambda_1-\lambda_2}\sum_{j=0}^k\binom{k}{j}\frac{\Gamma(k-j+\lambda_2-\lambda_1)}{\Gamma(\lambda_2-\lambda_1)}(-1)^js^{j}\partial^j\mu\geq 0.
\end{align*}
\hfill$\square$

\noindent
Stamatis Koumandos\\
Department of Mathematics and Statistics\\
The University of Cyprus \\ 
P.\ O.\ Box 20537\\
1678 Nicosia, Cyprus\\
{\em email}:\hspace{2mm}{\tt skoumand@ucy.ac.cy}

\vspace{0.5cm}

\noindent
Henrik Laurberg Pedersen\\
Department of Mathematical Sciences\\
University of Copenhagen \\
Universitetsparken 5\\
DK-2100, Denmark\\
{\em email}:\hspace{2mm}{\tt henrikp@math.ku.dk}

\end{document}